\documentclass[12pt]{article}
\usepackage{color}
\definecolor{darkblue}{rgb}{0.00,0.25,0.50}
\usepackage[colorlinks,filecolor=blue,citecolor=darkblue]{hyperref}

\usepackage{amsfonts,amssymb,amsmath,amsthm}
\usepackage{url}
\usepackage{enumerate}
\usepackage[ukrainian, russian, english]{babel}
\usepackage[cp1251]{inputenc}

\usepackage{cite}
\begin{document}

\selectlanguage{ukrainian} \thispagestyle{empty}

\title{}

\begin{center}
\textbf{\Large НЕРІВНОСТІ  ТИПУ ЛЕБЕГА ДЛЯ СУМ ВАЛЛЕ ПУССЕНА ТА ЇХ
ІНТЕРПОЛЯЦІЙНИХ АНАЛОГІВ НА КЛАСАХ
$(\psi,\bar{\beta})$-ДИФЕРЕНЦІЙОВНИХ ФУНКЦІЙ}
\end{center}
\vskip0.5cm
\begin{center}
В.\,А.~Войтович, А.\,П.~Мусієнко\\ \emph{\small Інститут
математики НАН України, Київ}
\end{center}
\vskip0.5cm

\begin{abstract}
We obtain the estimates of steady rates of deviations of the de
Vall\'{e}e Poussin sums and interpolation analogues of sums of
Vall\'{e}e Poussin from the functions that belong to the space
$C_{\bar{\beta}}^\psi L_s, \ 1\leq s\leq\infty$ and are
represented through the best approximations of
$(\psi,\bar{\beta})$-differentiable functions of this sort by
trigonometric polynomials in the metric $L_s$.
 \vskip 1.5mm

Одержано оцінки  норм відхилень сум Валле Пуссена та їх
інтерполяційних аналогів  від функцій з множин
$C_{\bar{\beta}}^\psi L_s, \ 1\leq s\leq\infty$, які виражаються
через найкращі наближення $(\psi,\bar{\beta})$-похідних таких
функцій тригонометричними поліномами в метриці простору $L_s$.
\end{abstract}

\vskip1cm

 Робота є продовженням досліджень \cite{Serd_2004,Serduk_Ovsii_2008,Serd_2009,Serd_2010,
Serduk_Ovsii_2011(1),Serd_Ovsii_2011,Serduk_Mysienko_2010,
Serduk_Mysienko_2012_1,Serduk_Mysienko_2012_2} по вивченню
апроксимативних властивостей сум Валле Пуссена або їх
інтерполяційних аналогів
\cite{Stepanec_Serdyuk_2000_i,Serdyuk_2002_i,Serdyuk_2005_c,Serdyuk_2012}
на класах $(\psi,\bar{\beta})$-диференційовних функцій.

Нехай $L_s$, $1\leq s<\infty$, --- простір $2\pi $-періодичних
сумовних в \mbox{$s$--му} степені на $(0,2\pi)$ функцій $f(t)$ з
нормою \mbox{$\|f\|_{s}=\Big(\int\limits_{-\pi
}^{\pi}{|f(t)|^s}\,dt\Big)^{1/s}.$} $L_\infty$ --- простір
вимірних і істотно обмежених $2\pi $-періодичних функцій $f(t)$ з
нормою \mbox{$\|f\|_{\infty}=\mathop{\rm ess\,sup}\limits_{t\
}|f(t)|.$} $C$ --- простір неперервних $2\pi$-періодичних функцій
$f(t),$ в якому норма задається рівністю
\mbox{$\|f\|_{C}=\max\limits_{t}{|f(t)|}.$}

Нехай $f$ --- $2\pi$-періодична, сумовна  функція (${f\in L_1}$) і
$$
\frac{a_0}{2}+\sum\limits_{k=1}^\infty(a_k\cos kx+b_k\sin
kx)
$$
--- її ряд Фур'є. Нехай, далі $\psi=\psi(k)$ і ${\bar{\beta}}=\beta_k, \
k\in\mathbb{N}$ --- довільні послідовності дійсних чисел. Якщо ряд
$$
\sum\limits_{k=1}^\infty\frac{1}{\psi(k)}\bigg(a_k\cos\big(kx+\frac{\beta_k\pi}{2}\big)+
b_k\sin\big(kx+\frac{\beta_k\pi}{2}\big)\bigg)
$$
є рядом Фур'є деякої сумовної функції $\varphi$, то цю функцію
називають $(\psi,\bar{\beta})$-похідною функції $f$ і позначають
через $f_{\bar{\beta}}^\psi$  \cite[c. 33]{Step_monog_1987}.
Множину всіх функцій $f$, які мають $(\psi,\bar{\beta})$--похідну,
позначають через $L_{\bar{\beta}}^\psi$. Якщо $f\in
L_{\bar{\beta}}^\psi$ і в той же час
$f_{\bar{\beta}}^\psi\in\mathfrak{N}$, де $\mathfrak{N}$
--- деяка підмножина з ${L_1^0=\{\varphi\in
L_1:\int\limits_{-\pi}^{\pi }\varphi(t) dt=0\}}$, то записують
$f\in L_{\bar{\beta}}^\psi\mathfrak{N}$. Якщо
$F_{\bar{\beta}}^\psi=f$, то функцію $F$ називають
$(\psi,{\bar{\beta}})$-інтегралом функції $f$, при цьому записують
$F(x)={\cal J}_{\bar{\beta}}^\psi(f;x)$. Покладемо
$C_{\bar{\beta}}^\psi=L_{\bar{\beta}}^\psi\bigcap C$,
$C_{\bar{\beta}}^\psi\mathfrak{N}=L_{\bar{\beta}}^\psi\mathfrak{N}\bigcap
C$. Означення $(\psi,\bar{\beta})$--похідних та
$(\psi,\bar{\beta})$--інтегралів та класів
$L^{\psi}_{\bar{\beta}}\mathfrak{N}$ та
$C^{\psi}_{\bar{\beta}}\mathfrak{N}$ належать О.І. Степанцю
\cite[c. 112]{Step_monog_2002(1)}.

Надалі будемо вважати, що послідовність $\psi(k)$, яка породжує
класи $L_{\bar{\beta}}^\psi\mathfrak{N}$ і
$C_{\bar{\beta}}^\psi\mathfrak{N}$,   задовольняє умову ${\cal
D}_0$ \ ($\psi\in{\cal D}_0$), тобто  $\psi(k)$ додатна і така, що
\begin{equation}\label{S4.0.4}
\mathop {\rm
\lim}\limits_{k\rightarrow\infty}\frac{\psi(k+1)}{\psi(k)}=0.
\end{equation}
 У випадку, коли $\psi\in {\cal D}_0,$ множини
$C^{\psi}_{\bar{\beta}}\mathfrak{N}$ складаються з функцій,
регулярних в усій комплексній площині, тобто з цілих функцій.

Відомо \cite[c. 144]{Step_monog_2002(1)}, що класи
$L_{\bar{\beta}}^\psi\mathfrak{N}$ складаються з функцій, які
майже при всіх $x\in\mathbb{R}$ можна зобразити у вигляді згортки
\begin{equation}\label{0.3}
f(x)=\frac{a_0}{2}+\frac{1}{\pi}\int\limits_{-\pi}^{\pi}\varphi(x-t)\Psi_{\bar{\beta}}(t)dt,
\ \varphi\in\mathfrak{N}, \ \varphi\perp1
\end{equation}
з сумовним ядром $\Psi_{\bar{\beta}}(t)$, ряд Фур'є  якого має
вигляд
$$
{\Psi_{\bar{\beta}}(t)\sim\sum\limits_{k=1}^\infty\psi(k)\cos\big(kt-\frac{\beta_k\pi}{2}\big), \ \psi(k)>0, \ \beta_k\in\mathbb{R}, \ k\in\mathbb{N}.}
$$
Якщо ж $f\in C_{\bar{\beta}}^\psi\mathfrak{N}$, то рівність
\eqref{0.3} виконується для всіх $x\in \mathbb{R}$.

\pagestyle{empty} \makeatletter
\renewcommand{\@evenhead}{\thepage \hfill {\it В.\,А.~Войтович,  А.\,П.~Мусієнко }}
\renewcommand{\@oddhead}{{\it {{{{Нерівності типу Лебега на класах}}}} \dots \hfill {\rm \thepage}}}
 \makeatother

Одиничну кулю  простору $L^0_s, \ 1\leq s\leq\infty,$ позначимо
через $U^0_s$ і покладемо $C_{\bar{\beta}}^\psi
U^0_s=C_{{\bar{\beta}},s}^\psi$, \ $L_{\bar{\beta}}^\psi
U^0_s=L_{{\bar{\beta}},s}^\psi$.

Нехай $\mathcal{T}_{2m-1}$ підпростір тригонометричних поліномів ${t_{m-1}(x)=\sum\limits_{k=0}^{m-1}(\alpha_k\cos kx+\gamma_k\sin
kx), \ \alpha_k, \ \gamma_k\in\mathbb{R}}$, порядок яких не
перевищує $m-1$. Величина
$$
E_m(f)_{X}=\mathop {\rm
\inf}\limits_{t_{m-1}\in{\mathcal{T}}_{2m-1}} \|f-t_{m-1}\|_X
$$
є найкращим наближенням функції $f\in X\subset L_1$ в метриці
простору $X$ тригонометричними поліномами порядку $m-1$. Далі в
ролі $X$ виступатимуть простори $C$ або $L_s, \  1\leq
s\leq\infty$.

Позначимо через  $V_{n,p}(f)$ суми Валле Пуссена 
 функції $f\in L_1$, тобто поліноми вигляду
$$
V_{n,p}(f)=V_{n,p}(f;x)=\frac{1}{p}\sum_{k=n-p}^{n-1}S_k(f;x),
$$
де $S_k(f)=S_k(f;x)$ --- частинні суми Фур'є  порядку $k$ функції
$f$, а $p=p(n)$ --- певний натуральний параметр, $p\leqslant n$.
При $p=1$ суми Валле Пуссена  $V_{n,p}(f)$ є  частинними сумами
Фур'є $S_{n-1}(f)$ порядку $n-1$, якщо ж $p=n$, то суми
$V_{n,p}(f)$ перетворюються у відомі суми Фейєра $\sigma_{n-1}(f)$
порядку $n-1$:
$$
\sigma_{n-1}(f)=\sigma_{n-1}(f;x)=\frac{1}{n}\sum_{k=0}^{n-1}S_k(f;x).
$$

Крім  звичайних сум Валле Пуссена $V_{n,p}(f)$ будемо розглядати
їх інтерполяційні аналоги $\widetilde{V}_{n,p}(f)$.

 Нехай $f \in C$. Через $\widetilde{S}_{n-1}(f;x)$ будемо позначати тригонометричний поліном порядку
$n-1$, що інтерполює $f(x)$ у точках
\mbox{$x^{(n-1)}_k=\frac{2k\pi}{2n-1}$,} $k \in \mathbb{Z}$, тобто
такий, що
$$
{\widetilde{S}}_{n-1}(f;x^{(n-1)}_k)=f(x^{(n-1)}_k),~ k \in \mathbb{Z}.
$$
Інтерполяційний тригонометричний поліном ${\widetilde{S}}_{n-1}(f;x)$,
 можна записати в такий спосіб (див., наприклад, \cite[c. 13--14]{zigmund}):
\begin{equation}
\label{sn}
 {
\widetilde{S}}_{n-1}(f;x)=\frac{a^{(n-1)}_0}{2}+\sum\limits^{n-1}_{k=1}(a^{(n-1)}_k\cos
kx +b^{(n-1)}_k\sin kx),
\end{equation}
 де
\begin{equation}
\label{akn}
 a^{(n-1)}_k=\frac{2}{2n-1}\sum\limits^{2n-2}_{j=0}
f(x^{(n-1)}_j) \cos kx^{(n-1)}_j ,~ k=0, 1,...,n,
\end{equation}
\begin{equation}
\label{bkn}
 b^{(n-1)}_k=\frac{2}{2n-1}\sum\limits^{2n-2}_{j=0}
f(x^{(n-1)}_j) \sin kx^{(n-1)}_j,~ k=1,2,...,n.
\end{equation}

Поліноми
\begin{equation}
\label{vp} {\mathop{V}\limits^{\sim}}_{n,p}(f;x)=
\frac{a^{(n-1)}_0}{2}\lambda^{(n)}_0+\sum\limits^{n-1}_{k=1}\lambda^{(n)}_k(a^{(n-1)}_k\cos
kx +b^{(n-1)}_k\sin kx),
\end{equation}
де
\begin{equation}
\label{lambda}
 \lambda_{k}^{(n)}=
  \begin{cases}
    1, & 0\leqslant k\leqslant n-p, \\
    1-\frac{k-n+p}{p}, &n-p+1\leqslant k\leqslant n.
  \end{cases}
\end{equation}
$p\in \mathbb{N}$, $1\leq{p}\leq{n}$, а $a^{(n-1)}_k$ і
$b^{(n-1)}_k$ означені формулами \eqref{akn} та \eqref{bkn}
відповідно, називають інтерполяційними аналогами сум Валле Пуссена
з параметрами $n$ та $p$. При $p=1$ суми
${\mathop{V}\limits^{\sim}}_{n,p}(f;x)$ співпадають з
інтерполяційними тригонометричними поліномами
${\mathop{S}\limits^{\sim}}_{n-1}(f;x)$. У випадку $p=n $  суми
${\mathop{V}\limits^{\sim}}_{n,p}(f;x)$ перетворюються в
інтерполяційні суми Фейєра $\widetilde{\sigma}_{n-1}(f;x)$ порядку
$n-1$:
$$
\widetilde{\sigma}_{n-1}(f;x)=\frac{1}{n}\sum\limits^{n-1}_{k=0}{\mathop{S}\limits^{\sim}}^{(n-1)}_k(f;x),
$$
де
$$
\widetilde{S}^{(n-1)}_k(f;x)=\frac{a^{(n-1)}_0}{2}+\sum\limits^{k}_{j=1}(a^{(n-1)}_j\cos
jx +b^{(n-1)}_j\sin jx).
$$
В загальному випадку інтерполяційні суми
$\widetilde{V}_{n,p}(f;x)$ виражаються через суми
$\widetilde{S}_k^{(n-1)}$ наступним чином
$$
\widetilde{V}_{n,p}(f;x)=\frac{1}{p}\sum\limits^{n-1}_{k=n-p}\widetilde{S}^{(n-1)}_k(f;x).
$$

Позначимо через $\rho_{n,p}(f;x)$ і $\widetilde{\rho}_{n,p}(f;x)$ величини вигляду
$$
\rho_{n,p}(f;x)=f(x)-V_{n,p}(f;x),
$$
$$
\widetilde{\rho}_{n,p}(f;x)=f(x)-\widetilde{V}_{n,p}(f;x).
$$

В роботі встановлено  асимптотично непокращувані аналоги
нерівностей типу Лебега для відхилень сум  $V_{n,p}(f)$ та
$\widetilde{V}_{n,p}(f;x)$ на множинах $C_{\bar{\beta}}^\psi L_s$
при $\psi\in {\cal D}_0$, $\beta_k\in\mathbb{R}, \ k\in\mathbb{N}$
і ${1\leq s\leq\infty}$.

\bf Теорема 1. \it{ Нехай $\psi\in {\cal D}_0$,
$\beta_k\in\mathbb{R}, \ k\in\mathbb{N}$. Тоді для довільних $f\in
C_{\bar{\beta}}^\psi L_s$, ${1\leq s<\infty}$, $n,p\in\mathbb{N}$,
$p\leq n$ справедлива нерівність
\begin{equation}\label{S4.T.2.1}
\|\rho_{n,p}(f;x)\|_C\leq$$$$\leq\Big(\frac{\|\cos
t\|_{s'}}{\pi p}\psi
(n-p+1)+O(1)\sum\limits_{k=n-p+2}^\infty\psi(k)\tau_{n,p}(k)\Big)E_{n-p+1}(f_{\bar{\beta}}^\psi)_{L_s}.
\end{equation}

При цьому для будь-якої функції $f\in C^\psi_{\bar{\beta}} L_s$,
${1\leq s<\infty}$ і довільних $n,p\in \mathbb{N}$, $p\leq n$ в
множині $C^\psi_{\bar{\beta}} L_s$,  ${1\leq s<\infty}$,
знайдеться функція ${F(x)=F(f;n;p;x)}$ така, що
${E_{n-p+1}{(F_{\bar{\beta}}^\psi)}_{L_s}=E_{n-p+1}{(f_{\bar{\beta}}^\psi)}_{L_s}}$,
і для неї при $n-p\rightarrow\infty$ виконується рівність
\begin{equation}\label{S4.T.2.2}
\|\rho_{n,p}(F;x)\|_C=$$$$=\Big(\frac{\|\cos t\|_{s'}}{\pi p}\psi
(n-p+1)+O(1)\sum\limits_{k=n-p+2}^\infty\psi(k)\tau_{n,p}(k)\Big)E_{n-p+1}(F_{\bar{\beta}}^\psi)_{L_s}.
\end{equation}
 У \eqref{S4.T.2.1} і \eqref{S4.T.2.2} $s'=\frac{s}{s-1}$, коефіцієнти
$\tau_{n,p}(k)$ визначаються рівністю
\begin{equation}\label{S4.T.1.9}
\tau_{n,p}(k)= \left\{\begin{array}{ll}
1-\frac{n-k}{p},& \ \ n-p+1\leq k\leq n-1,\\
1,&  \ \  k\geq n,
\end{array}\right.
\end{equation}
а $O(1)$ --- величини, рівномірно обмежені відносно всіх
розглядуваних параметрів.

{\textbf{\textit{Доведення теореми {\rm\bf1}.}}} \ \rm Доведення
будемо проводити за схемою, яка запропонована в теоремі 1 роботи
\cite{Serduk_Mysienko_2012_2}. Нехай $f\in C_{\bar{\beta}}^\psi
L_s$, ${1\leq s\leq\infty}$. Тоді в кожній точці $x\in\mathbb{R}$
(див., наприклад, \cite[c. 810]{Rykasov_2003}) має місце
інтегральне зображення
\begin{equation}\label{S4.T.2.4}
\rho_{n,p}(f;x)=\frac{1}{\pi}\int\limits_{-\pi}^\pi
f_{\bar{\beta}}^\psi(x-t)\Psi_{1,n,p}(t)dt=$$$$=\frac{\psi(n-p+1)}{\pi
p}\int\limits_{-\pi }^{\pi
}f_{\bar{\beta}}^\psi(x-t)\cos\Big((n-p+1)t-\frac{\beta_{n-p+1}\pi}{2}\Big)dt+$$$$+
\frac{1}{\pi}\int\limits_{-\pi }^{\pi
}f_{\bar{\beta}}^\psi(x-t)\Psi_{2,n,p}(t)dt,
\end{equation}
де $\Psi_{j,n,p}(t)$ означається рівністю
\begin{equation}\label{S4.T.1.8}
\Psi_{j,n,p}(t)=\sum\limits_{k=n-p+j}^\infty\tau_{n,p}(k)\psi(k)\cos\Big(kt-\frac{\beta_k\pi}{2}\Big),
\ j\in\mathbb{N}.
\end{equation}

Функції $\cos\big((n-p+1)t-\frac{\beta_{n-p+1}\pi}{2}\big)$ та
$\Psi_{2,n,p}(t)$ ортогональні до будь-якого тригонометричного
полінома $t_{n-p}$ порядку не вищого  ${n-p}$, тому в силу
$\eqref{S4.T.2.4}$
\begin{equation}\label{S4.T.2.5}
\rho_{n,p}(f;x)=$$$$=\frac{\psi(n-p+1)}{\pi p}\int\limits_{-\pi
}^{\pi
}\delta_{n,p}(x-t)\cos\Big((n-p+1)t-\frac{\beta_{n-p+1}\pi}{2}\Big)dt+$$$$+
\frac{1}{\pi}\int\limits_{-\pi }^{\pi
}\delta_{n,p}(x-t)\Psi_{2,n,p}(t)dt,
\end{equation}
де
\begin{equation}\label{S4.T.1.16}
\delta_{n,p}(\cdot)=f_{\bar{\beta}}^\psi(\cdot)-t_{n-p}(\cdot).
\end{equation}
Обравши в \eqref{S4.T.2.5} у ролі $t_{n-p}(\cdot)$ поліном
$t_{n-p}^\ast(\cdot)$ найкращого наближення в просторі
 $L_s$ функції $f_{\bar{\beta}}^\psi(\cdot)$ та використовуючи формулу \eqref{S4.T.1.9} і
 нерівність Гельдера
\begin{equation}\label{S4.T.1.18'}
\bigg\|\int\limits_{-\pi}^\pi K(t-u)\varphi(u)du\bigg\|_C\leq
$$$$\leq\|K\|_{s'}\|\varphi\|_s, \ \varphi\in L_s, \ K\in L_{s'},
\ 1\leq s\leq\infty, \ \frac{1}{s}+\frac{1}{s'}=1,
\end{equation}
(див., наприклад, \cite[c. 43]{Korn}),  інтеграли рівності
\eqref{S4.T.2.5} оцінимо наступним чином:
\begin{equation}\label{S4.T.2.9}
\big\|\int\limits_{-\pi }^{\pi }\delta_{n,p}(x-t)\cos(n-p+1)t\,
dt\big\|_C\leq$$$$\leq\|\delta_{n,p}\|_s\|\cos t\|_{s'}\leq\|\cos
t\|_{s'}E_{n-p+1}(f_{\bar{\beta}}^\psi)_{L_s},
\end{equation}
\begin{equation}\label{S4.T.2.7}
\big\|\frac{1}{\pi}\int\limits_{-\pi }^{\pi
}\delta_{n,p}(x-t)\Psi_{2,n,p}(t)dt\big\|_C\leq$$$$\leq
\frac{1}{\pi}\|\delta_{n,p}\|_s\|\Psi_{2,n,p}\|_{s'}\leq
\frac{2^{1/s'}}{\pi^{1/s}}\sum\limits_{k=n-p+2}^{\infty}\psi(k)\tau_{n,p}(k)E_{n-p+1}(f_{\bar{\beta}}^\psi)_{L_s}.
\end{equation}
Об'єднуючи \eqref{S4.T.2.9} і \eqref{S4.T.2.7} отримуємо
\eqref{S4.T.2.1}.

Доведемо другу частину теореми. З інтегрального зображення
\eqref{S4.T.2.4} і факту ортогональності функції $\Psi_{2,n,p}(t)$
до будь-якого тригонометричного полінома
$t_{n-p}\in\mathcal{T}_{2(n-p)+1}$ випливає, що для довільної
функції $f\in C_{\bar{\beta}}^\psi L_s, \ \ 1\leq s<\infty$,
$\psi\in{\cal D}_0$, $\beta_k\in\mathbb{R}$, $k\in\mathbb{N}$
виконується рівність
\begin{equation}\label{S4.T.2.10}
|\rho_{n,p}(f;x)|=$$$$=\frac{\psi(n-p+1)}{\pi
p}\Big|\int\limits_{-\pi }^{\pi
}f_{\bar{\beta}}^\psi(x-t)\cos\big((n-p+1)t-\frac{\beta_{n-p+1}\pi}{2}\big)dt\Big|+$$$$+
O(1)\sum\limits_{k=n-p+2}^\infty\psi(k)\tau_{n,p}(k)E_{n-p+1}(f_{\bar{\beta}}^\psi)_{L_s}.
\end{equation}
Враховуючи \eqref{S4.T.2.10}, щоб переконатися в справедливості
\eqref{S4.T.2.2} досить показати, що якою б не була функція
$\varphi\in L_s, \ \ 1\leq s<\infty$ знайдеться функція
${\Phi(\cdot)=\Phi(\varphi;\cdot)}$, для якої при всіх
$n,p\in\mathbb{N}$, $p\leq n$
\begin{equation}\label{S4.T.2.11}
E_{n-p+1}(\Phi)_{L_s}=E_{n-p+1}(\varphi)_{L_s}
\end{equation}
і, крім того, має місце рівність
\begin{equation}\label{S4.T.2.12}
\big|\int\limits_{-\pi }^{\pi
}\Phi(t)\cos\Big((n-p+1)t+\frac{\beta_{n-p+1}\pi}{2}\big)dt\Big|=\|\cos
t\|_{s'}E_{n-p+1}(\varphi)_{L_s}.
\end{equation}

В якості $\Phi(\cdot)$  розглянемо функцію
\begin{equation}\label{S4.T.2.13}
\Phi(t)=\|\cos
t\|_{s'}^{1-s'}\big|\cos\big((n-p+1)t+\frac{\beta_{n-p+1}\pi}{2}\big)\big|^{s'-1}\times$$$$\times
\mathrm{sign}\cos\big((n-p+1)t+\frac{\beta_{n-p+1}\pi}{2}\big)E_{n-p+1}(\varphi)_{L_s}.
\end{equation}
Для неї
\begin{equation}\label{S4.T.2.14}
\|\Phi(t)\|_s=\|\cos t\|_{s'}^{1-s'}\Big(\int\limits_{-\pi }^{\pi
}\big|\cos\big((n-p+1)t+$$$$+\frac{\beta_{n-p+1}\pi}{2}\Big)\big|^{(s'-1)s}dt\big)^\frac{1}{s}E_{n-p+1}(\varphi)_{L_s}=$$$$
=\|\cos
t\|_{s'}^{1-s'}\big\|\cos\big((n-p+1)t-\frac{\beta_{n-p+1}\pi}{2}\big)\big\|^{s'-1}_{s'}E_{n-p+1}(\varphi)_{L_s}=$$$$=E_{n-p+1}(\varphi)_{L_s}.
\end{equation}
Крім того, оскільки для довільного
$t_{n-p}\in\mathcal{T}_{2(n-p)+1}$
$$
\int\limits_{-\pi}^{\pi}t_{n-p}(\tau)|\Phi(\tau)|^{s-1}\mathrm{sign}\Phi(\tau)d\tau=
\Big(\|\cos
t\|_{s'}^{1-s'}E_{n-p+1}(\varphi)_{L_s}\Big)^{s-1}\times$$$$\times\int\limits_{-\pi}^{\pi}t_{n-p}(\tau)\cos\Big((n-p+1)\tau-\frac{\beta_{n-p+1}\pi}{2}\Big)d\tau=0,
$$
то на підставі теореми 1.4.5 роботи \cite[c. 28]{Korn}  робимо
висновок, що поліном ${t_{n-p}^\ast\equiv0}$ є поліномом
найкращого наближення функції $\Phi(t)$ в метриці простору ${L_s,
1\leq s<\infty}$. Отже,  з урахуванням \eqref{S4.T.2.14},
\begin{equation}\label{S4.T.2.14'}
E_{n-p+1}(\Phi)_{L_s}=\|\Phi\|_s=E_{n-p+1}(\varphi)_{L_s}.
\end{equation}
В силу \eqref{S4.T.2.13},
$$
\Big|\int\limits_{-\pi}^{\pi}\Phi(t)\cos\big((n-p+1)t+\frac{\beta_{n-p+1}\pi}{2}\big)dt\Big|=$$$$=
\|\cos
t\|_{s'}^{1-s'}E_{n-p+1}(\varphi)_{L_s}\Big|\int\limits_{-\pi}^{\pi}\big|\cos\big((n-p+1)t+\frac{\beta_{n-p+1}\pi}{2}\big)\big|^{s'-1}
\times$$$$\times\mathrm{sign}\cos\big((n-p+1)t+\frac{\beta_{n-p+1}\pi}{2}\big)\cos\big((n-p+1)t+\frac{\beta_{n-p+1}\pi}{2}\big)dt\Big|=$$$$=
\|\cos
t\|_{s'}^{1-s'}E_{n-p+1}(\varphi)_{L_s}\int\limits_{-\pi}^{\pi}\big|\cos\big((n-p+1)t+\frac{\beta_{n-p+1}\pi}{2}\big)\big|^{s'}dt=$$$$=\|\cos
t\|_{s'}E_{n-p+1}(\varphi)_{L_s}.
$$
З останніх  співвідношень випливає \eqref{S4.T.2.12}. Теорему 1
доведено.

\vskip 3mm

\bf Теорема 2. \it{ Нехай $\psi\in {\cal D}_0$,
$\beta_k\in\mathbb{R}$, $k\in\mathbb{N}$. Тоді для довільних $f\in
C_{\bar{\beta}}^\psi L_\infty$, $n,p\in\mathbb{N}$, $p\leq n$
справедлива нерівність
\begin{equation}\label{S4.T.3.1(1)}
\|\rho_{n,p}(f;x)\|_C\leq\frac{1}{p}\bigg(\frac{4}{\pi}\psi(n-p+1)+O(1)\Big(\frac{\psi^2(n-p+2)}{\psi(n-p+1)}
+$$$$+p\sum\limits_{k=n-p+3}^\infty\psi(k)\tau_{n,p}(k)\Big)\bigg)E_{n-p+1}(f_{\bar{\beta}}^\psi)_{L_\infty}.
\end{equation}

При цьому для будь-якої функції $f\in C_{\bar{\beta}}^\psi
L_\infty$ і довільних ${n,p\in \mathbb{N}}$, $p\leq n$  знайдеться
функція ${{F(x)=F(f;n;p;x)}\in C_{\bar{\beta}}^\psi C}$ така, що
${E_{n-p+1}{(F_{\bar{\beta}}^\psi)}_C=E_{n-p+1}{(f_{\bar{\beta}}^\psi)}_{L_\infty}}$
і для неї при ${n-p\rightarrow\infty}$ виконується рівність
\begin{equation}\label{S4.T.3.2(1)}
\|\rho_{n,p}(F;x)\|_C=\frac{1}{p}\bigg(\frac{4}{\pi}\psi(n-p+1)+O(1)\Big(\frac{\psi^2(n-p+2)}{\psi(n-p+1)}
+$$$$+p\!\!\!\!\sum\limits_{k=n-p+3}^\infty\psi(k)\tau_{n,p}(k)\Big)\bigg)E_{n-p+1}(F_{\bar{\beta}}^\psi)_C.
\end{equation}
У \eqref{S4.T.3.1(1)} і \eqref{S4.T.3.2(1)} коефіцієнти
$\tau_{n,p}(k)$ визначаються рівністю \eqref{S4.T.1.9}, а
$O(1)$~--- величини, рівномірно обмежені відносно всіх
розглядуваних параметрів.

{\textbf{\textit{Доведення теореми {\rm\bf2}}}} \ \rm  будемо
проводити за схемою, яка запропонована в теоремі 2 роботи
\cite{Serduk_Mysienko_2012_2}. Нехай $f\in C_{\bar{\beta}}^\psi
L_\infty,$ $\psi\in{\cal D}_0$. Виходячи з \eqref{S4.T.2.4}, та
враховуючи факт ортогональності функції $\Psi_{1,n,p}(t)$ до
будь-якого полінома $t_{n-p}$ порядку не вищого за $n-p$, можемо
записати
\begin{equation}\label{S4.T.3.3}
\rho_{n,p}(f;x)=\frac{1}{\pi}\int\limits_{-\pi}^\pi
\delta_{n,p}(x-t)\Psi_{1,n,p}(t)dt=$$$$=\frac{1}{\pi
p}\int\limits_{-\pi}^\pi\delta_{n,p}(x-t)\big(\psi(n-p+1)
\cos\big((n-p+1)t-\frac{\beta_{n-p+1}\pi}{2}\big)+$$$$+2\psi(n-p+2)\cos\big((n-p+2)t-\frac{\beta_{n-p+2}\pi}{2}\big)+$$$$+
p\!\!\!\!\sum\limits_{k=n-p+3}^\infty\!\!\!\tau_{n,p}(k)\psi(k)\cos\big(kt-\frac{\beta_k\pi}{2}\big)\big)dt,
\end{equation}
де $\tau_{n,p}(k)$ і $\delta_{n,p}(\cdot)$ визначаються рівностями
\eqref{S4.T.1.9} і \eqref{S4.T.1.16} відповідно.

Обравши в \eqref{S4.T.3.3} у ролі $t_{n-p}$ поліном $t_{n-p}^\ast$
найкращого наближення в просторі $L_\infty$ функції
$f_{\bar{\beta}}^\psi$ і застосувавши нерівність
\eqref{S4.T.1.18'} при $s=\infty$, отримуємо оцінку
\begin{equation}\label{S4.T.3.7}
\|\rho_{n,p}(f;x)\|_C\leq\frac{1}{\pi p}\big(\big\|\psi(n-p+1)
\cos\big((n-p+1)t-\frac{\beta_{n-p+1}\pi}{2}\big)+$$$$+2\psi(n-p+2)\cos\big((n-p+2)t-\frac{\beta_{n-p+2}\pi}{2}\big)\big\|_1+$$$$+p\big\|\sum\limits_{k=n-p+3}^\infty\tau_{n,p}(k)\psi(k)\cos\big(kt-\frac{\beta_k\pi}{2}\big)\big\|_1\big)E_{n-p+1}{(f_{\bar{\beta}}^\psi)}_{L_\infty}\leq$$$$
\leq\frac{1}{\pi p}\Big(\big\|\psi(n-p+1)
\cos(n-p+1)t+$$$$+2\psi(n-p+2)\cos\big((n-p+2)t+\alpha_{\bar{\beta},n,p}\big)\big\|_1+$$$$
+O(1)p\!\!\!\!\sum\limits_{k=n-p+3}^\infty\tau_{n,p}(k)\psi(k)\Big)E_{n-p+1}{(f_{\bar{\beta}}^\psi)}_{L_\infty},
\end{equation}
де
\begin{equation}\label{lamda}
\alpha_{\bar{\beta},n,p}=\frac{\beta_{n-p+2}\pi}{2}-
    \frac{n-p+2}{n-p+1}\frac{\beta_{n-p+1}\pi}{2}.
\end{equation}
Як випливає з роботи С.О.~Теляковського \cite[c.
512--513]{Teliacovs_1989},
\begin{equation}\label{S4.T.3.8}
\big\|\psi(n-p+1)\cos(n-p+1)t+2\psi(n-p+2)\cos\big((n-p+2)t+\alpha_{\bar{\beta},n,p}\big)\big\|_1+$$$$+O(1)p\!\!\!\!\sum\limits_{k=n-p+3}^\infty\tau_{n,p}(k)\psi(k)\leq$$$$\leq
4\psi(n-p+1)+O(1)\Big(\frac{\psi^2(n-p+2)}{\psi(n-p+1)}+p\!\!\!\!\sum\limits_{k=n-p+3}^\infty\!\!\!\tau_{n,p}(k)\psi(k)\Big).
\end{equation}
Співвідношення \eqref{S4.T.3.7} і \eqref{S4.T.3.8} доводять
нерівність \eqref{S4.T.3.1(1)}.

Доведемо другу частину теореми. Виходячи з інтегрального
зображення \eqref{S4.T.3.3} і використовуючи факт ортогональності
функції
${\sum\limits_{k=n-p+3}^\infty\!\!\!\tau_{n,p}(k)\psi(k)\cos\big(kt-\frac{\beta_k\pi}{2}\big)}$
до будь-якого тригонометричного полінома $t_{n-p}$ порядку не
вищого $n-p$, для довільної функції $f$  з множини
$C_{\bar{\beta}}^\psi L_\infty$, $\psi\in {\cal D}_0$,
$\beta_k\in\mathbb{R}, \ k\in\mathbb{N}$ виконується рівність
\begin{equation}\label{S4.T.3.9}
|\rho_{n,p}(f;x)|=$$$$=\frac{\psi(n-p+1)}{\pi
p}\Big|\int\limits_{-\pi}^\pi f_{\bar{\beta}}^\psi(x-t)\big(
\cos\big((n-p+1)t-\frac{\beta_{n-p+1}\pi}{2}\big)+$$$$+2\frac{\psi(n-p+2)}{\psi(n-p+1)}\cos\big((n-p+2)t-\frac{\beta_{n-p+2}\pi}{2}\big)\big)dt\Big|+$$$$+
O(1)\!\!\!\sum\limits_{k=n-p+3}^\infty\!\!\!\tau_{n,p}(k)\psi(k)E_{n-p+1}{(f_{\bar{\beta}}^\psi)}_{L_\infty}.
\end{equation}
Для доведення \eqref{S4.T.3.2(1)}, з урахуванням \eqref{S4.T.3.9},
досить встановити, що для довільної ${\varphi\in
L_\infty^0=\{\varphi\in L_\infty: \varphi\bot1\}}$ існує функція
${\Phi(\cdot)=\Phi(\varphi;\cdot)\in C}$ для якої при всіх $n,p\in
\mathbb{N}$, $p\leq n$
$$
E_{n-p+1}(\Phi)_C=E_{n-p+1}(\varphi)_{L_\infty}
$$
і, крім того, при $n-p\rightarrow\infty$ має місце рівність
\begin{equation}\label{S4.T.3.10}
\big|\int\limits_{-\pi}^\pi \Phi(t)\big(
\cos\big((n-p+1)t+\frac{\beta_{n-p+1}\pi}{2}\big)+$$$$+2\frac{\psi(n-p+2)}{\psi(n-p+1)}\cos\big((n-p+2)t+\frac{\beta_{n-p+2}\pi}{2}\big)\big)dt\big|=$$$$=
\bigg(4+O(1)\Big(\frac{\psi(n-p+2)}{\psi(n-p+1)}\Big)^2\bigg)E_{n-p+1}(\varphi)_{L_\infty}.
\end{equation}

Покладемо
$$
\varphi_0(t)=\mathrm{sign}\cos\Big((n-p+1)t+\frac{\beta_{n-p+1}\pi}{2}\Big)E_{n-p+1}(\varphi)_{L_\infty}
$$
і через $\varphi_\delta(t)$ позначимо $2\pi$-періодичну функцію,
яка збігається з $\varphi_0(t)$ скрізь, за виключенням
$\delta$-околів (${0<\delta<\frac{\pi}{2(n-p+1)}}$) точок
${t_k=\frac{(2k+1-\beta_{n-p+1})\pi}{2(n-p+1)}, \ \
k\in\mathbb{Z}}$, де вона лінійна і її графік сполучає точки
$(t_k-\delta, \varphi_0(t_k-\delta))$ і $(t_k+\delta,
\varphi_0(t_k+\delta))$. Функція $\varphi_\delta(t)$ неперервна і
у точках ${\tau_k=\frac{(2k-\beta_{n-p+1})\pi}{2(n-p+1)}, \ \
k=1,2,...,2(n-p+1)}$ періоду
${\big(-\frac{\beta_{n-p+1}\pi}{2(n-p+1)},
2\pi-\frac{\beta_{n-p+1}\pi}{2(n-p+1)}\big]}$ досягає по
абсолютній величині максимального значення, яке дорівнює
$E_{n-p+1}(\varphi)_{L_\infty}$,
 почергово змінюючи знак. Тому її поліном найкращого рівномірного наближення порядку не вищого $n-p$,
 згідно з критерієм Чебишова, є поліном, що тотожно дорівнює нулю і, отже,
\begin{equation}\label{S4.T.3.11}
E_{n-p+1}(\varphi_\delta)_C=\|\varphi_\delta\|_C=E_{n-p+1}(\varphi)_{L_\infty}.
\end{equation}
Враховуючи \eqref{S4.T.1.18'} і \eqref{S4.T.3.11}, одержуємо
\begin{equation}\label{S4.T.3.12}
\big|\int\limits_{-\pi}^\pi \varphi_\delta(t)\big(
\cos\big((n-p+1)t+\frac{\beta_{n-p+1}\pi}{2}\big)+$$$$+2\frac{\psi(n-p+2)}{\psi(n-p+1)}\cos\big((n-p+2)t+\frac{\beta_{n-p+2}\pi}{2}\big)\big)dt\big|\leq$$$$\leq
\int\limits_{-\pi}^\pi
\big|\cos(n-p+1)t+2\frac{\psi(n-p+2)}{\psi(n-p+1)}\cos\big((n-p+2)t+\alpha_{\bar{\beta},n,p}\big)\big|dt
\times$$$$\times E_{n-p+1}(\varphi)_{L_\infty},
\end{equation}
де $\alpha_{\bar{\beta},n,p}$ визначається рівністю \eqref{lamda}.
Із нерівності (19) роботи \cite{Teliacovs_1989}, випливає оцінка
\begin{equation}\label{S4.T.3.13}
\int\limits_{-\pi}^\pi
\big|\cos(n-p+1)t+2\frac{\psi(n-p+2)}{\psi(n-p+1)}\cos\big((n-p+2)t+\alpha_{\bar{\beta},n,p}\big)\big|dt\leq$$$$\leq
4+O(1)\Big(\frac{\psi(n-p+2)}{\psi(n-p+1)}\Big)^2.
\end{equation}
З іншого боку,
\begin{equation}\label{S4.T.3.14}
\big|\int\limits_{-\pi}^\pi \varphi_\delta(t)\big(
\cos\big((n-p+1)t+\frac{\beta_{n-p+1}\pi}{2}\big)+$$$$+2\frac{\psi(n-p+2)}{\psi(n-p+1)}\cos\big((n-p+2)t+\frac{\beta_{n-p+2}\pi}{2}\big)\big)dt\big|=$$$$
=\big|\int\limits_{-\pi}^\pi \varphi_0(t)\big(
\cos\big((n-p+1)t+\frac{\beta_{n-p+1}\pi}{2}\big)+$$$$+2\frac{\psi(n-p+2)}{\psi(n-p+1)}\cos\big((n-p+2)t+\frac{\beta_{n-p+2}\pi}{2}\big)\big)dt\big|+O(1)r_{n,p}(\delta),
\end{equation}
де
\begin{equation}\label{S4.T.3.15}
r_{n,p}(\delta)=\big|\int\limits_{-\pi}^\pi
(\varphi_\delta(t)-\varphi_0(t))\big(
\cos\big((n-p+1)t+\frac{\beta_{n-p+1}\pi}{2}\big)+$$$$+2\frac{\psi(n-p+2)}{\psi(n-p+1)}\cos\big((n-p+2)t+\frac{\beta_{n-p+2}\pi}{2}\big)\big)dt\big|.
\end{equation}
Оскільки $\psi\in {\cal D}_0$, то для досить великих номерів $n-p$
справджується нерівність  $\frac{\psi(n-p+2)}{\psi(n-p+1)}<1$ і,
отже,
\begin{equation}\label{S4.T.3.16}
r_{n,p}(\delta)<3\int\limits_{-\pi}^{\pi}|\varphi_\delta(t)-\varphi_0(t)|dt\leq6(n-p+1)\delta
E_{n-p+1}(\varphi)_{L_\infty},
\end{equation}
то вибравши $\delta$ настільки малим, щоб виконувалась умова
\begin{equation}\label{S4.T.3.17}
0<\delta<\frac{1}{n-p+1}\Big(\frac{\psi(n-p+2)}{\psi(n-p+1)}\Big)^2,
\end{equation}
із \eqref{S4.T.3.16} одержимо оцінку
\begin{equation}\label{S4.T.3.18}
r_{n,p}(\delta)=O(1)\Big(\frac{\psi(n-p+2)}{\psi(n-p+1)}\Big)^2E_{n-p+1}(\varphi)_{L_\infty}.
\end{equation}
Оскільки
$\int\limits_{-\pi}^{\pi}\varphi_0(t)\cos\big((n-p+2)t+\frac{\beta_{n-p+2}\pi}{2}\big)dt=0$,
то врахувавши розклад функцій $\mathrm{sign}\cos(n-p+1)t$ та
$\mathrm{sign}\sin(n-p+1)t$ в ряд Фур'є отримаємо
\begin{equation}\label{S4.T.3.19}
\big|\int\limits_{-\pi}^\pi \varphi_0(t)\big(
\cos\big((n-p+1)t+\frac{\beta_{n-p+1}\pi}{2}\big)+$$$$+2\frac{\psi(n-p+2)}{\psi(n-p+1)}\cos\big((n-p+2)t+\frac{\beta_{n-p+2}\pi}{2}\big)\big)dt\big|=$$$$
=\big|\int\limits_{-\pi}^{\pi}\varphi_0(t)\cos\big((n-p+1)t+\frac{\beta_{n-p+1}\pi}{2}\big)dt\big|=$$$$=
\int\limits_{-\pi}^{\pi}\big|\cos\big((n-p+1)t+\frac{\beta_{n-p+1}\pi}{2}\big)\big|dt
E_{n-p+1}(\varphi)_{L_\infty}=4E_{n-p+1}(\varphi)_{L_\infty}.
\end{equation}
Із формул \eqref{S4.T.3.12}--\eqref{S4.T.3.14}, \eqref{S4.T.3.18}
і \eqref{S4.T.3.19} випливає, що для функції
${\Phi(t)=\varphi_\delta(t)}$ у якій параметр $\delta$ задовольняє
умову \eqref{S4.T.3.17}, при $n-p+1\rightarrow\infty$ має місце
рівність \eqref{S4.T.3.10}, а отже, і \eqref{S4.T.3.2(1)}. Теорему
2 доведено.

Далі, розглянемо  аналоги теорем 1 та 2 для величини
$\widetilde{\rho}_{n,p}(f;x)$.

У випадку $p=1$, тобто коли суми $\widetilde{V}_{n,p}(f;x)$ є
інтерполяційними поліномами $\widetilde{S}_{n-1}(f;x)$, нерівності
типу Лебега на класах цілих функцій встановлені в роботі
\cite{Serdyuk_2012}. Тому ми розглянемо лише випадок $2\leq p\leq
n$.

\bf Теорема 3. \it{ Нехай $\psi\in {\cal D}_0$,
$\beta_k\in\mathbb{R}, \ k\in\mathbb{N}$. Тоді для довільних $f\in
C_{\bar{\beta}}^\psi L_s$, ${1\leq s<\infty}$, $n,p\in\mathbb{N}$,
$2\leq p\leq n$ справедлива нерівність
\begin{equation}\label{TV1}
|\widetilde{\rho}_{n,p}(f;x)|\leq$$$$\leq\Big(\frac{\|\cos
t\|_{s'}}{\pi p}\psi
(n-p+1)+O(1)\sum\limits_{k=n-p+2}^\infty\psi(k)\tau_{n,p}(k)\Big)E_{n-p+1}(f_{\bar{\beta}}^\psi)_{L_s}.
\end{equation}
При цьому для будь-якої функції $f\in C^\psi_{\bar{\beta}} L_s$,
${1\leq s<\infty}$ і довільних $n,p\in \mathbb{N}$, $p\leq n$ в
множині $C^\psi_{\bar{\beta}} L_s$,  ${1\leq s<\infty}$,
знайдеться функція ${F(x)=F(f;n;p;x)}$ така, що
${E_{n-p+1}{(F_{\bar{\beta}}^\psi)}_{L_s}=E_{n-p+1}{(f_{\bar{\beta}}^\psi)}_{L_s}}$,
і для неї при $n-p\rightarrow\infty$ виконується рівність
\begin{equation}\label{TV1a}
|\widetilde{\rho}_{n,p}(F;x)|=$$$$=\Big(\frac{\|\cos t\|_{s'}}{\pi
p}\psi
(n-p+1)+O(1)\sum\limits_{k=n-p+2}^\infty\psi(k)\tau_{n,p}(k)\Big)E_{n-p+1}(F_{\bar{\beta}}^\psi)_{L_s}.
\end{equation}
У \eqref{TV1} та \eqref{TV1a} $s'=\frac{s}{s-1}$, коефіцієнти
$\tau_{n,p}(k)$ означаються рівністю \eqref{S4.T.1.9}, а
$O(1)$~--- величини, рівномірно обмежені відносно всіх
розглядуваних параметрів.}

\bf Теорема 4. \it{ Нехай $\psi\in {\cal D}_0$,
$\beta_k\in\mathbb{R}, \ k\in\mathbb{N}$. Тоді для довільних $f\in
C_{\bar{\beta}}^\psi L_{\infty}$ і будь--яких $n,p\in\mathbb{N}$,
$2\leq p\leq n$ справедлива нерівність
\begin{equation}\label{TV2}
|\widetilde{\rho}_{n,p}(f;x)|\leq$$$$\leq\Big(\frac{4}{\pi p}\psi
(n-p+1)+O(1)\sum\limits_{k=n-p+2}^\infty\psi(k)\tau_{n,p}(k)\Big)E_{n-p+1}(f_{\bar{\beta}}^\psi)_{L_{\infty}}.
\end{equation}
При цьому для будь-якої функції $f\in C^\psi_{\bar{\beta}}
L_{\infty}$ і довільних $n,p\in \mathbb{N}$, $p\leq n$ в множині
$C^\psi_{\bar{\beta}}C$, знайдеться функція ${F(x)=F(f;n;p;x)}$
така, що
${E_{n-p+1}{(F_{\bar{\beta}}^\psi)}_{C}=E_{n-p+1}{(f_{\bar{\beta}}^\psi)}_{L_{\infty}}}$,
і для неї при $n-p\rightarrow\infty$ виконується рівність
\begin{equation}\label{TV2a}
|\widetilde{\rho}_{n,p}(F;x)|=$$$$=\Big(\frac{4}{\pi p}\psi
(n-p+1)+O(1)\sum\limits_{k=n-p+2}^\infty\psi(k)\tau_{n,p}(k)\Big)E_{n-p+1}(F_{\bar{\beta}}^\psi)_{C},
\end{equation}
де коефіцієнти $\tau_{n,p}(k)$ означаються рівністю
\eqref{S4.T.1.9}, а $O(1)$~--- величини, рівномірно обмежені
відносно всіх розглядуваних параметрів.}

\rm

Оскільки доведення теорем 3 та 4 не відрізняються, тому проведемо
їх  разом.

{\textbf{\textit{Доведення теорем {\rm\bf3} та {\rm\bf4}.}}} \rm
Нехай $f\in C^{\psi}_{\bar{\beta}}L_s$, ${1\leq s\leq\infty}$,
${\psi\in{\cal D}_0}$. В лемі 2 роботи \cite{Vojtovich_serdyuk}
встановлено, що коли $\psi(k)>0$,
$\sum\limits^\infty_{k=1}\psi(k)<\infty,~\beta_k \in \mathbb{R}, \
k\in\mathbb{N}$, то для довільної функції $f \in
C^\psi_{\bar{\beta}}L_s, \ 1\leq s \leq\infty$, в кожній точці $x
\in \mathbb{R}$ мають місце рівності
\begin{equation}
\label{1l2}
\widetilde{\rho}_{n,p}(f;x)=\rho_{n,p}(f;x)+O(1)E_{n-p+1}(f^\psi_{\bar{\beta}})_{L_s}\sum\limits^\infty_{k=n}\psi(k),
\end{equation}де $O(1)$ --- величина,
рівномірно обмежена відносно всіх розглядуваних параметрів.

З урахуванням формули \eqref{S4.T.2.5} запишемо \eqref{1l2} в
такому вигляді:
\begin{equation}
\label{DT3_1} \widetilde{\rho}_{n,p}(f;x)=$$$$=\frac{\psi(n-p+1)}{\pi
p}\int\limits_{-\pi }^{\pi
}\delta_{n,p}(x-t)\cos\Big((n-p+1)t-\frac{\beta_{n-p+1}\pi}{2}\Big)dt+$$$$+
\frac{1}{\pi}\int\limits_{-\pi }^{\pi
}\delta_{n,p}(x-t)\Psi_{2,n,p}(t)dt+O(1)E_{n-p+1}(f^\psi_{\bar{\beta}})_{L_s}\sum\limits^\infty_{k=n}\psi(k).
\end{equation}
Застосувавши нерівність
\begin{equation}
\label{nerKorn}
 \int\limits^{\pi}_{-\pi}\varphi(t)K(t)dt \leq
\|\varphi\|_s\|K\|_{s'},
\varphi \in L_s, \  K \in
L_{s'}, \frac{1}{s}+\frac{1}{s'}=1,~1\leq s\leq \infty,
\end{equation}
(див., наприклад, \cite[с. 391]{Korn}), та рівності
\eqref{S4.T.2.9} та \eqref{S4.T.2.7}, отримаємо  \eqref{TV1} та
\eqref{TV2}.

Доведемо тепер рівності \eqref{TV1a} та \eqref{TV2a}. Розглянемо
спочатку випадок $1\leq s <\infty$. Як показано в доведенні
теореми 1, для функції ${\Phi(\cdot)=\Phi(\varphi;\cdot)}$,
означеної рівністю \eqref{S4.T.2.13}, при всіх $n,p\in\mathbb{N}$,
$p\leq n$
$$
E_{n-p+1}(\Phi)_{L_s}=E_{n-p+1}(f^{\psi}_{\bar{\beta}})_{L_s}
$$
і, крім того, має місце рівність
\begin{equation}\label{DTV1}
|\rho_{n,p}(F;x)|=$$$$=\Big(\frac{\|\cos t\|_{s'}}{\pi p}\psi
(n-p+1)+O(1)\sum\limits_{k=n-p+2}^\infty\psi(k)\tau_{n,p}(k)\Big)E_{n-p+1}(\Phi)_{L_s},
\end{equation}
де $F={\cal J}^{\psi}_{\bar{\beta}}\Phi$. Тому з \eqref{1l2} та
\eqref{DTV1} випливає рівність \eqref{TV1a}.

Розглянемо випадок, коли $s=\infty$. Як показано в доведенні
теореми 2, для функції
${\varphi_{\delta}(\cdot)=\varphi_{\delta}(\varphi;\cdot)}$,
означеної в доведенні другої частини теореми 2, при всіх
$n,p\in\mathbb{N}$, $p\leq n$
$$
E_{n-p+1}(\Phi)_{L_s}=E_{n-p+1}(f^{\psi}_{\bar{\beta}})_{L_s}
$$
 і, крім того, має місце рівність
\begin{equation}\label{DTV2}
|\rho_{n,p}(F;x)|=\frac{1}{p}\bigg(\frac{4}{\pi}\psi(n-p+1)+O(1)\Big(\frac{\psi^2(n-p+2)}{\psi(n-p+1)}
+$$$$+p\!\!\!\!\sum\limits_{k=n-p+3}^\infty\psi(k)\tau_{n,p}(k)\Big)\bigg)E_{n-p+1}(\varphi_{\delta})_C,
\end{equation}
де $F={\cal J}^{\psi}_{\bar{\beta}}\varphi_{\delta}$. Тому з
\eqref{1l2} та \eqref{DTV2} випливає рівність \eqref{TV1a}.
 Теореми 3 та 4 доведено.

\rm Оскільки для сум
$\sum\limits_{k=n-p+j}^\infty\tau_{n,p}(k)\psi(k)$, які фігурують
в теоремах {\rm 1 --- 4,}  мають місце рівності
\begin{equation}\label{S4.T.3.24'}
   \sum\limits_{k=n-p+j}^{\infty}\tau_{n,p}(k)\psi(k)=$$$$=\left\{\begin{array}{ll}
   \sum\limits_{k=n-p+j}^{n-1}\frac{k-n+p}{p}\psi (k)+\sum\limits_{k=n}^{\infty}\psi (k), & \ p>j, \ \ j\in\mathbb{N},\\
   \sum\limits_{k=n-p+j}^{\infty}\psi(k), & \ p\leq j, \ \ j\in\mathbb{N},
  \end{array}\right.
\end{equation}
то, як неважко переконатися, для них справедлива наступна оцінка
зверху:
$$
\sum\limits_{k=n-p+j}^{\infty}\tau_{n,p}(k)\psi(k)\leqslant$$$$\leq
\min\big\{\sum\limits_{k=n-p+j}^{\infty}\!\!\!\psi
(k),\frac{1}{p}\sum\limits_{k=n-p+j}^{\infty}(k-n+p)\psi
(k)\big\},\ \ \ \ j\in \mathbb{N}.
$$
Отже, в співвідношеннях \eqref{S4.T.2.1} і \eqref{S4.T.2.2}
теореми 1, співвідношеннях \eqref{S4.T.3.1(1)} і
\eqref{S4.T.3.2(1)} теореми 2, співвідношеннях \eqref{TV1} і
\eqref{TV1a} теореми 3 та співвідношеннях \eqref{TV2} і
\eqref{TV2a} теореми 4 величини
$O(1)\!\!\!\sum\limits_{k=n-p+j}^\infty\tau_{n,p}(k)\psi(k)$ можна
замінити на
${O(1)\min\big\{\sum\limits_{k=n-p+j}^{\infty}\!\!\!\psi
(k),\frac{1}{p}\sum\limits_{k=n-p+j}^{\infty}(k-n+p)\psi
(k)\big\}}$,  ${j=2,3}$.

{\textbf{\textit{Зауваження.}}} При $\beta_k=\beta, \
k\in\mathbb{N}$  теореми 1 та 2 встановлені в роботі
\cite{Serduk_Mysienko_2012_2}].

\renewcommand{\refname}{}
\makeatletter\renewcommand{\@biblabel}[1]{#1.}\makeatother

\end{document}